\documentclass[11pt]{amsart}
\usepackage{amsfonts} 
\usepackage{amsmath} 
\usepackage{amsthm} 
\usepackage{amssymb} 
\usepackage{enumitem} 

\newcommand{\sep}{\itemsep 0pt}

\DeclareMathOperator{\rank}{rank}
\DeclareMathOperator{\soc}{soc}
\DeclareMathOperator{\End}{End}

\newtheorem{theorem}{Theorem}[section]
\newtheorem{proposition}[theorem]{Proposition}
\newtheorem{lemma}[theorem]{Lemma}
\newtheorem{corollary}[theorem]{Corollary}

\theoremstyle{definition}
\newtheorem{definition}[theorem]{Definition}
\newtheorem{example}[theorem]{Example}
\newtheorem{remark}[theorem]{Remark}

\begin{document}

\title[Rank of elements of general rings]%
{Rank of elements of general rings in connection with unit-regularity}

\author{Nik Stopar}
\address{Faculty of Electrical Engineering, University of Ljubljana, Tr\v za\v ska cesta 25, 1000 Ljubljana, Slovenia}
\email{nik.stopar@fe.uni-lj.si}

\begin{abstract}
We define the rank of elements of general unital rings, discuss its properties and give several examples to support the definition. In semiprime rings we give a characterization of rank in terms of invertible elements. As an application we prove that every element in the socle of a unital semiprime ring is unit-regular.
\end{abstract}

\maketitle

{\footnotesize \emph{Key Words:} rank, minimal right ideal, idempotent, unit-regular element

}

\section{Introduction}
Although rank is mostly associated with matrices and linear algebra, there are several other areas where the rank of an element has been defined in a more general setting. Two of the most important examples are:
\smallskip
\begin{itemize}[leftmargin=*]\itemsep 5pt
\item[-] The theory of polynomial identities (PI) and generalized polynomial identities (GPI), where the rank is defined in the setting of primitive rings. It gives an alternative way of describing the socle of a primitive ring which plays a fundamental role in this theory. For example, the condition of finite rank appears in several important theorems of Amitsur \cite{Ami}. For a full exposition on this subject see \cite{Bei-Mar-Mik,Row} and the references therein.
\item[-] The theory of Banach algebras, where rank was investigated mostly in the setting of semisimple Banach algebras. Perhaps the most thorough investigation of rank in this setting was done by Aupetit and Mouton \cite{Aup-Mou} and later by Bre\v sar and \v Semrl \cite{Bre-Sem}, although the definition itself appeared even earlier in several other papers (see \cite{Mou-Rau,Puh} and the references therein).
Aupetit and Mouton defined rank via the spectrum and used it as a tool to define the trace and determinant of certain elements of the algebra. They relied heavily on analytical methods. On the other hand, Bre\v sar and \v Semrl gave a completely algebraic definition of rank and presented several equivalent characterizations, showing in particular, that their definition coincides with those given earlier.
Later Brits, Lindeboom and Raubenheimer \cite{Bri-Lin-Rau} also used rank to describe certain properties of the Drazin inverse and the Drazin index of elements of Banach algebras.
\end{itemize}
\medskip

As far as we are aware, in the most general setting of arbitrary (unital) rings, the rank has not been investigated yet. The aim of this paper is to do precisely that. We believe that rank may be a useful tool for understanding the structure of rings, especially semiprime rings, as indicated by our results below.

Throughout this paper, $R$ will be a unital associative ring (or possibly a unital associative algebra over a field). Following \cite{Bre-Sem} we define the right rank of an element $a \in R$ as follows.

\begin{definition}\label{def_rank}
An element $a \in R$ has \emph{right rank} $0$ if and only if $a=0$. An element $a \in R$ has \emph{right rank} $1$ if and only if $a \neq 0$ and $a$ is contained in some minimal right ideal of $R$. An element $a \in R$ has \emph{right rank} $n>1$ if and only if $a$ is contained in a sum of $n$ minimal right ideals of $R$, but is not contained in any sum of less than $n$ minimal right ideals of $R$. An element $a \in R$ has \emph{infinite right rank} if and only if $a$ is not contained in any sum of minimal right ideals of $R$. The right rank of $a \in R$ will be denoted by $\rank_r a$.
\end{definition}

Similarly we define the \emph{left rank} of an element $a \in R$ and denote it by $\rank_l a$. With the convention that the sum of zero minimal right ideals is $0$, we can say, that the \emph{right rank} of an element $a \in R$ is the least nonnegative integer $n$, such that $a$ is contained in a sum of $n$ minimal right ideals of $R$. If such an integer does not exist, then the right rank of $a$ is infinite.

Our paper is organized as follows. In \S 2 we recall some definitions and give some notation that will be used throughout the paper. In \S 3 we present several interesting examples to justify and support the definition of rank given above.
Our key results about rank are contained in \S 4, where we describe decompositions of elements as sums of elements of rank one. We investigate in details the rank of idempotents and give some alternative characterizations of rank.
Most of these results are valid in general rings, however, our main result that characterizes the rank in terms of invertible elements (see Theorem~\ref{idem_rank} and Corollary~\ref{rank_char}) is most definitive in the setting of semiprime rings. As an application of our main theorem, we prove that the socle of a semiprime ring is unit-regular (see Theorem~\ref{unit_regular}).

Another application of the results obtained here will be presented in a subsequent paper \cite{Sto2}, where we will investigate the structure of generalized corner rings.

\section{Preliminaries} 
For a unital ring $R$ we denote by $M_n(R)$ the ring of all $n \times n$ matrices with entries in $R$. More generally, the set of all $n \times m$ matrices with entries in $R$ will be denoted by $M_{n,m}(R)$.
Standard matrix units in $M_n(R)$ will be denoted by $E_{ij}$, $1 \leq i,j \leq n$, so $E_{ij}$ is a matrix whose only nonzero entry is entry $(i,j)$ and is equal to $1$. We shall also often need the ring of all upper triangular $n \times n$ matrices over $R$, which we will denote by $T_n(R)$.

For the group of all multiplicatively invertible elements of a ring $R$ we will be using the standard notation $\mathcal{U}(R)$.
Recall that an element $a \in R$ is called \emph{regular} if there exists an element $b \in R$ such that $a=aba$. In this case $ab$ and $ba$ are idempotents.
If there exists $b \in \mathcal{U}(R)$ that satisfies this condition, then the element $a$ is called \emph{unit-regular}. Equivalently, $a$ is unit-regular if there exists $x \in \mathcal{U}(R)$ and an idempotent $e \in R$, such that $a=ex$.

The \emph{right socle} of a ring $R$, denoted by $\soc R_R$, is defined as the sum of all minimal right ideals of $R$ and is a two sided ideal of $R$. In particular, $\soc R_R=0$ if $R$ has no minimal right ideals.
By our Definition~\ref{def_rank}, $\soc R_R$ is precisely the set of all elements of $R$ of finite right rank. The \emph{left socle} of $R$, denoted by $\soc \phantom{}_R R$, is defined analogously via left ideals. If $R$ is a semiprime ring, then the left and the right socle of $R$ coincide and are thus simply called the \emph{socle} of $R$ and denoted by $\soc R$.

\section{Examples}

In this section we give some examples to illustrate and justify the definition of rank (see Definition~\ref{def_rank}). We start with matrices.

\begin{example}
Let $F$ be a field and consider the ring $R=M_n(F)$. This ring is endowed with the usual matrix rank defined as the dimension of the image of a matrix. It turns out that this usual matrix rank coincides with both left and right rank as defined in Definition~\ref{def_rank}. We shall verify this in the next more general example. This, in particular, justifies the terminology.
\end{example}

\begin{remark}\label{sum}
Observe that, since we are dealing with unital rings, an element $a \in R$ has right rank $1$ if and only if the right ideal $aR$ is minimal. An element of right rank $n>1$ can be expressed as a sum of $n$ elements of right rank $1$, but not as a sum of less than $n$ elements of right rank $1$. This property is in fact equivalent to the definition of right rank $n>1$.
\end{remark}

Next we consider matrices over arbitrary division rings.

\begin{example}\label{ex_matrix}
Let $D$ be a division ring and $R=M_n(D)$. In this case, one has to be a bit more careful when defining the matrix rank via dimension. Recall that the \emph{row rank} of a matrix $A \in R$ is defined as the dimension of the left vector space over $D$ spanned by the rows of $A$.
The \emph{column rank} of $A$ is defined as the dimension of the right vector space over $D$ spanned by the columns of $A$. It is well known (see \cite{Hun} for details) that the row and column rank of a matrix coincide and are thus simply called the \emph{rank} of a matrix. Again, this rank coincides with both left and right rank as defined in Definition~\ref{def_rank}.
To see this, it is enough to show, that any minimal right ideal of $R$ is generated by a matrix of column rank $1$, and any matrix of column rank $1$ generates a minimal right ideal of $R$. Let $K$ be a minimal right ideal of $R$ and choose $0 \neq A \in K$.
We may assume that $A$ is in its column echelon form, since this form is achieved by multiplying the matrix from the right by an appropriate invertible matrix and the result is still in $K$. Then $AE_{11}R$ is a nonzero right ideal of $R$ contained in $K$.
Hence $K=AE_{11}R$ and $K$ is generated by a matrix $AE_{11}$ of column rank $1$. Now take a matrix $B$ of column rank $1$. Again we may assume that $B$ is in its column echelon form, since this does not change the right ideal generated by $B$.
Since the column rank of $B$ is $1$, this means that $B=BE_{11}$. Choose any nonzero matrix $C \in BR$ and write it as $C=BX=BE_{11}X$ for some matrix $X \in R$. Then $E_{11}X \neq 0$. Since only the first row of $E_{11}X$ is nonzero, its reduced column echelon form is $E_{11}$.
Thus there exists an invertible matrix $U$ such that $E_{11}XU=E_{11}$. This implies $CU=B$ and hence $CR=BR$. This shows that $BR$ is a minimal right ideal of $R$.
\end{example}

Recall that the left and right socle in a general ring need not coincide. This indicates that the left and right rank of an element need not be equal. In fact, the following example demonstrates that there is completely no connection between the two, apart from the trivial connection in rank $0$.

\begin{example}\label{mn_rank}
Let $m$ and $n$ be two positive integers and $F$ a field. Let $R$ be the subalgebra of $M_{2mn}(F)$ consisting of all the matrices of the form
\[\left[\begin{array}{cccc|cccc}
A & 0 & \cdots & 0 & B_{11} & B_{12} & \cdots & B_{1m} \\
0 & A & \cdots & 0 & B_{21} & B_{22} & \cdots & B_{2m} \\
\vdots & \vdots & \ddots & \vdots & \vdots & \vdots & \ddots & \vdots \\
0 & 0 & \ldots & A & B_{n1} & B_{n2} & \cdots & B_{nm} \\
\hline
0 & 0 & \cdots & 0 & C & 0 & \cdots & 0 \\
0 & 0 & \cdots & 0 & 0 & C & \cdots & 0 \\
\vdots & \vdots & \ddots & \vdots & \vdots & \vdots & \ddots & \vdots \\
0 & 0 & \cdots & 0 & 0 & 0 & \cdots & C \\
\end{array}\right]=
\left[\begin{array}{c|c}
\mathcal{A} & \mathcal{B} \\
\hline
0 & \mathcal{C} \\
\end{array}\right],\]
where $A$ is an arbitrary $m \times m$ matrix, $C$ is an arbitrary $n \times n$ matrix, and $B_{ij}$ are arbitrary $m \times n$ matrices. So $\mathcal{A}$, $\mathcal{B}$, and $\mathcal{C}$ are $mn \times mn$ matrices, $\mathcal{B}$ is arbitrary, while $\mathcal{A}$ and $\mathcal{C}$ are block diagonal with constant blocks of size $m \times m$ and $n \times n$ respectively. Take the matrices
\[J=\left[\begin{array}{c|c} 0 & \mathcal{I} \\ \hline 0 & 0 \\ \end{array}\right],\quad K=\left[\begin{array}{c|c} \mathcal{I} & 0 \\ \hline 0 & 0 \\ \end{array}\right],\quad L=\left[\begin{array}{c|c} 0 & 0 \\ \hline 0 & \mathcal{I} \\ \end{array}\right] \quad\in R \]
where $\mathcal{I}$ is the $mn \times mn$ identity matrix. We claim that
\[\begin{array}{lcl}
\rank_r J=n, && \rank_l J=m \\
\rank_r K=\infty, & &\rank_l K=m, \\
\rank_r L=n, && \rank_l L=\infty.
\end{array}\]
Due to the symmetry, it is enough to consider the right rank. Observe that the set
\[\Delta=\{\left[\begin{array}{c|c} 0 & \mathcal{B} \\ \hline 0 & \mathcal{C} \\ \end{array}\right] \in R \}\]
is a right ideal of $R$ and $S \Delta \neq 0$ for any nonzero $S \in R$. This implies that any minimal right ideal of $R$ is contained in $\Delta$. Hence, $\rank_r K=\infty$. Now take arbitrary
\[\left[\begin{array}{c|c} \mathcal{A} & \mathcal{B} \\ \hline 0 & \mathcal{C} \\ \end{array}\right] \in R \qquad \textup{and} \qquad \left[\begin{array}{c|c} 0 & \mathcal{U} \\ \hline 0 & \mathcal{V} \\ \end{array}\right] \in \Delta,\]
where $\mathcal{A}$, $\mathcal{B}$, and $\mathcal{C}$ are as above, and similarly $\mathcal{U}=[U_{ij}]_{i,j}$, and $\mathcal{V}$ is block diagonal with blocks $V$. Then we have
\begin{equation}\label{right-action}
\left[\begin{array}{c|c} 0 & \mathcal{U} \\ \hline 0 & \mathcal{V} \\ \end{array}\right]\cdot \left[\begin{array}{c|c} \mathcal{A} & \mathcal{B} \\ \hline 0 & \mathcal{C} \\ \end{array}\right]=\left[\begin{array}{c|c} 0 & \mathcal{UC} \\ \hline 0 & \mathcal{VC} \\ \end{array}\right].
\end{equation}
Observe that multiplying $\mathcal{U}$ from the right by $\mathcal{C}$ means multiplying every block-column of $\mathcal{U}$ from the right by $C$. Multiplying $\mathcal{V}$ by $\mathcal{C}$ means multiplying $V$ by $C$. Now break the matrix $\mathcal{U}$ into its block-columns (of width $n$), put these block-columns one below the other to form a $m^2n \times n$ matrix, then attach matrix $V$ to the bottom of this matrix to get a $(m^2n+n) \times n$ matrix, and denote this matrix by $\mathcal{U} \boxplus \mathcal{V}$.
Then by the above observations, the right action \eqref{right-action} of $R$ on $\Delta$ is equivalent to the right action $(\mathcal{U} \boxplus \mathcal{V}) \cdot C$ of $M_n(F)$ on $M_{m^2n+n,n}(F)$.
Hence the right ideals of $R$ contained in $\Delta$ are in a one-to-one correspondence with the submodules of the right $M_n(F)$-module $M_{m^2n+n,n}(F)$. It is not hard to see, that the minimal submodules of $M_{m^2n+n,n}(F)$ are those, that are generated by a rank $1$ matrix (cf. Example~\ref{ex_matrix}). Therefore, the right rank of a matrix
\[\left[\begin{array}{c|c} 0 & \mathcal{U} \\ \hline 0 & \mathcal{V} \\ \end{array}\right] \in \Delta\]
is just the usual matrix rank of the matrix $\mathcal{U} \boxplus \mathcal{V}$. It is now clear that $\rank_r J=n$ and $\rank_r L=n$. It also follows from the above, that
\[\soc R_R=\{\left[\begin{array}{c|c}
0 & \mathcal{B} \\
\hline
0 & \mathcal{C} \\
\end{array}\right] \in R\}
\qquad \textup{and} \qquad
\soc \phantom{}_RR=\{\left[\begin{array}{c|c}
\mathcal{A} & \mathcal{B} \\
\hline
0 & 0 \\
\end{array}\right] \in R\}.\]
\end{example}

Before giving further examples, we observe that the rank satisfies the expected rank inequalities.

\begin{proposition}\label{rank_ineq}
For every $a,b \in R$ we have:
\begin{enumerate}\sep
\item $\rank_r(a+b) \leq \rank_r a+\rank_r b$,
\item $\rank_r(ab) \leq \min\{\rank_r a,\rank_r b\}$.
\end{enumerate}
\end{proposition}

\proof The first part follows directly from the definition.

For the second part we first show that if either of $a$ or $b$ has right rank $1$, then $ab$ has right rank $\leq 1$. If $a$ has right rank $1$, then $aR$ is a minimal right ideal. So either $ab=0$ or $ab \in aR$ is again an element of right rank $1$.
If $b$ has right rank $1$, then $bR$ is a minimal right ideal. If $ab \neq 0$, then $abR$ is again a minimal right ideal. This is because the map $bR \to abR$, $x \mapsto ax$, is a bijective right $R$-module homomorphism ($bR$ being a minimal right ideal, the kernel is either $0$ or $bR$, and it is not $bR$ since $ab \neq 0$).

Now assume $\rank_r a=n$, where $0<n<\infty$. Then by Remark~\ref{sum}, $a=a_1+a_2+\ldots+a_n$ for some $a_i$ with right rank $1$. Hence, $ab=a_1b+a_2b+\ldots+a_nb$, and by the above, all $a_ib$ have right rank $\leq 1$. The first part of proposition now implies that $\rank_r ab \leq \rank_r a$.
This last inequality holds also if the right rank of $a$ is $0$ or $\infty$. The same argument shows that $\rank_r ab \leq \rank_r b$. \endproof

\begin{corollary}\label{rank_inv}
For every $a \in R$ and $u \in \mathcal{U}(R)$ we have
\[\rank_r(au)=\rank_r(ua)=\rank_r a.\]
\end{corollary}

We will be using Proposition~\ref{rank_ineq} and Corollary~\ref{rank_inv} throughout the paper, often without explicit reference.

Example~\ref{mn_rank} shows that any combination of positive left and right rank is possible, so the two are completely independent in general. The situation is much simpler if $R$ is a semiprime ring, in which case the left and right rank coincide. This is a consequence of the following well known fact from the theory of idempotents, the proof of which can be found in \cite[Corollary~10.23, Proposition~21.16]{Lam}.

\begin{proposition}\label{min_idem}
Let $R$ be a ring and $e \in R$ an idempotent.
\begin{enumerate}
\item If $eR$ (resp. $Re$) is a minimal right (resp. left) ideal of $R$, then $eRe$ is a division ring. The converse holds if $R$ is a semiprime ring.
\item If $R$ is a semiprime ring, then every minimal right (resp. left) ideal of $R$ is of the form $eR$ (resp. $Re$) for some idempotent $e \in R$.
\end{enumerate}
\end{proposition}

Thus, in a semiprime ring $R$, every minimal one-sided ideal is generated by an idempotent, and $eR$ is a minimal right ideal if and only if $Re$ is a minimal left ideal. This, together with Proposition~\ref{rank_ineq}, easily implies that right rank $1$ and left rank $1$ coincide, and by Remark~\ref{sum} the left and right rank in $R$ coincide.

\begin{corollary}\label{semiprime}
If $R$ is a semiprime ring, then the left and right rank in $R$ coincide.
\end{corollary}

Idempotents that generate minimal right ideals are usually called \emph{right irreducible} idempotents (see for example Lam \cite{Lam}). In our context, these are just the idempotents of right rank $1$.

Observe that by \cite[Remark~3.2]{Tug} every element in the socle of a semiprime ring is regular. In a subsequent paper \cite{Sto2} we prove a generalization of Corollary~\ref{semiprime}, namely that the left and right rank of a regular element of any ring coincide as long as they are both finite.

Next example considers rank in prime rings.

\begin{example}
Let $R$ be a prime ring with nonzero socle and $L$ a minimal left ideal in $R$. Then $R$ is a primitive ring with faithful simple left $R$-module $L$ (see \cite[Theorem~11.11]{Lam}). By Schur's Lemma, $D=\End_R L$ is a division ring and $L$ is a left vector space over $D$.
By the Jacobson Density Theorem, there is an inclusion $\varphi: R \hookrightarrow \End_D L$ whose image is a dense subring of $\End_D L$. By \cite[Lemma~5.29]{Bre}, an element $a \in R$ has rank $1$ if and only if the operator $\varphi(a)$ has operator rank $1$.
Hence the rank of any element $r \in R$ as defined in Definition~\ref{def_rank} coincides with the operator rank of $\varphi(r) \in \End_D L$. In other words (cf. \cite[\S 7.1]{Row}), the rank of an element $r \in R$ in this case is the dimension of $rL$ as a left vector space over $D$. The same is true if we replace $L$ by any faithful simple left $R$-module.
\end{example}

\begin{example}
Let $A$ be a semisimple Banach algebra. There are several equivalent characterizations of rank known in this case, some involve spectral conditions and some involve representations of Banach algebras. For details we refer the reader to \cite{Aup-Mou,Bre-Sem} and the references therein. In \cite{Bre-Sem} it is shown that all these characterizations are equivalent to Definition~\ref{def_rank}.
\end{example}

\section{Properties of rank}

Let $R$ be a unital ring. The property described in Remark~\ref{sum} will often be used, so it is convenient to make the following definition.

\begin{definition}\label{mrd}
We say that $a=a_1+a_2+\ldots+a_n$ is a \emph{minimal right decomposition} of $a$, if all $a_i$ have right rank $1$ and $n$ is the right rank of $a$.
\end{definition}

It should come as no surprise that the rank is most well behaved on idempotent elements. So we first investigate minimal right decompositions of idempotents. We will need the following technical lemma.

\begin{lemma}\label{zero_div}
If $a=a_1+a_2+\ldots+a_n$ is a minimal right decomposition of $a$ and $ab=0$, then $a_ib=0$ for all $i$.
\end{lemma}

\proof Suppose $ab=0$ but say $a_1b \neq 0$. Then $a_1b=-(a_2+\ldots+a_n)b$. Since $a_1$ has right rank $1$, $a_1R$ is a minimal right ideal. Hence, $a_1b \neq 0$ implies $a_1bR=a_1R$, so there exists $x \in R$ such that $a_1bx=a_1$. But then $a_1=-(a_2+\ldots+a_n)bx$, so that $a=a_1+a_2+\ldots+a_n=(a_2+\ldots+a_n)(1-bx)$. By Proposition~\ref{rank_ineq}, this implies that the right rank of $a$ is at most $n-1$, which is a contradiction. \endproof

We remark, that the following argument from the proof of Lemma~\ref{zero_div}, will be often used without further explanation. If $a$ has right rank $1$ and $ab \neq 0$, then there exists an element $x$, such that $abx=a$.

We say that a set $\{a_1,a_2,\ldots,a_n\}$ of elements of $R$ is an \emph{orthogonal system}, if $a_ia_j=0$ for all $i \neq j$.

\begin{proposition}\label{idem_decomp}
If $e=a_1+a_2+\ldots+a_n$ is a minimal right decomposition of an idempotent $e$, then $\{a_1,a_2,\ldots,a_n\}$ is an orthogonal system of idempotents. 
\end{proposition}

\proof Observe that $e=e^2=ea_1+ea_2+\ldots+ea_n$. Suppose that $a_1(1-a_1)\neq 0$. Then there exists $y \in R$ such that $a_1(1-a_1)y=a_1$, which implies $a_1=(1-a_1)a_1y$. Hence $e=(e-ea_1)+ea_1=(e-ea_1)+e(1-a_1)a_1y=(e-ea_1)(1+a_1y)=(ea_2+\ldots+ea_n)(1+a_1y)$. This leads to a contradiction $\rank_r e \leq n-1$. Thus $a_1(1-a_1)=0$ and $a_1$ is an idempotent.
Clearly, $e=ea_1+ea_2+\ldots+ea_n$ is again a minimal right decomposition of $e$, which means that $e-ea_1=ea_2+\ldots+ea_n$ must be a minimal right decomposition of $e-ea_1$. Since $a_1$ is an idempotent, we have $(e-ea_1)a_1=e(1-a_1)a_1=0$. Hence, Lemma~\ref{zero_div} implies $ea_ia_1=0$ for all $i \neq 1$.
By symmetry, all $a_i$ are idempotents and $ea_ia_j=0$ for all $i \neq j$. By Lemma~\ref{zero_div}, $a_i(1-e)=0$ for all $i$, so that $a_ie=a_i$ for all $i$. Therefore $a_ia_j=a_i^2a_j=(a_ie)a_ia_j=a_i(ea_ia_j)=0$ for all $i \neq j$. We conclude that $\{a_1,a_2,\ldots,a_n\}$ is an orthogonal system of idempotents. \endproof

The converse of Proposition~\ref{idem_decomp} also holds.

\begin{corollary}\label{idem_decomp_2}
If $\{e_1,e_2,\ldots,e_n\}$ is an orthogonal system of idempotents of right rank $1$, then $e_1+e_2+\ldots+e_n$ is an idempotent of right rank $n$.
\end{corollary}

This is a direct consequence of the following more general proposition. Observe that for an element $a$ of right rank $1$, either $a^2=0$ or $a$ is not nilpotent.

\begin{proposition}\label{orthog_decomp}
If $\{e_1,e_2,\ldots,e_n\}$ is an orthogonal system of non-\-nil\-po\-tent elements of right rank $1$, then $e_1+e_2+\ldots+e_n$ is a non-nilpotent element of right rank $n$.
\end{proposition}

\proof We prove this by induction on $n$. For $n=1$ the conclusion is obvious. Suppose the right rank of $e=e_1+e_2+\ldots+e_n$ is $k<n$ and let $e=a_1+a_2+\ldots+a_k$ be a minimal right decomposition of $e$. Since $e_1^2=ee_1=a_1e_1+a_2e_1+\ldots+a_ke_1$, there exists $j$ such that $a_je_1 \neq 0$. By renumbering $a_i$'s, we may assume $j=1$, so that $a_1e_1 \neq 0$. This implies $a_1e_1x=a_1$ for some $x \in R$. Hence $a_1=a_1e_1x-(e-e_1)e_1x=e_1^2x-(e-a_1)e_1x$. Denote $f=e_2+e_3+\ldots+e_n$. Since $e_1f=0$, we have 
\begin{eqnarray*}
f^2 &=& ef=a_1f+(e-a_1)f=(e_1^2x-(e-a_1)e_1x)f+(e-a_1)f= \\
&=& e_1^2xf+(e-a_1)(1-e_1x)f=e_1^2xf+(a_2+a_3+\ldots+a_k)(1-e_1x)f.
\end{eqnarray*}
Now denote $b_1=e_1^2xf$ and $b_j=a_j(1-e_1x)f$ for $j \geq 2$. Then $f^2=b_1+b_2+\ldots+b_k$ and all $b_i$ have right rank at most $1$. Observe that $fe_1=0$ implies $fb_1=0$ and $b_jb_1=0$ for all $j$. Thus
\begin{eqnarray*}
f^4 &=& f^2(b_1+b_2+\ldots+b_k)=f^2(b_2+b_3+\ldots+b_k)= \\
&=& b_1(b_2+b_3+\ldots+b_k)+(b_2+b_3+\ldots+b_k)^2.
\end{eqnarray*}
On the other hand
\[b_1f^2=b_1(b_1+b_2+\ldots+b_k)=b_1(b_2+\ldots+b_k).\]
Subtracting the last two equalities we get $f^4-b_1f^2=(b_2+b_3+\ldots+b_k)^2$, which implies that the right rank of $f^4-b_1f^2$ is at most $k-1\leq n-2$. Now
\begin{eqnarray*}
f^4-b_1f^2 &=& (e_2^4+e_3^4+\ldots+e_n^4)-b_1(e_2^2+e_3^2+\ldots+e_n^2)= \\
&=& (e_2^2-b_1)e_2^2+(e_3^2-b_1)e_3^2+\ldots+(e_n^2-b_1)e_n^2= \\
&=& c_2+c_3+\ldots+c_n,
\end{eqnarray*}
where $c_i=(e_i^2-b_1)e_i^2$ for all $i \geq 2$. By induction, $f=e_2+e_3+\ldots+e_n$ has rank $n-1$, so this is its minimal right decomposition. In view of $fb_1=0$, Lemma~\ref{zero_div} implies $e_ib_1=0$ for all $i \geq 2$. Consequently, $e_i(e_j^2-b_1)=0$ for all $i,j \geq 2$, $i \neq j$, which implies $c_ic_j=0$ for all $i,j \geq 2$, $i \neq j$.
In addition, $e_ic_i^2=e_i(e_i^2-b_1)e_i^2(e_i^2-b_1)e_i^2=e_i^9$ for all $i \geq 2$. In particular, the elements $c_i$ are not nilpotent and have right rank $1$. By induction, the element $f^4-b_1f^2=c_2+c_3+\ldots+c_n$ has right rank $n-1$. This is a contradiction because we have already seen that this element has right rank at most $n-2$. We conclude that the right rank of $e$ is $n$.

For a positive integer $m$ we have $e^m=e_1^m+e_2^m+\ldots+e_n^m$. Since $\{e_1^m, e_2^m, \ldots, e_n^m\}$ is again an orthogonal system of non-nilpotent elements of right rank $1$, we infer from the above that the right rank of $e^m$ is $n$. In particular, $e^m \neq 0$. \endproof

The assumption that the elements are non-nilpotent in Proposition~\ref{orthog_decomp} is essential. For example, in the ring $M_n(\mathbb{C})$, the matrix $E_{1n}+E_{2n}+\ldots+E_{n-1,n}$ has rank $1$, although the summands form an orthogonal system.

In view of the fact, that Proposition~\ref{orthog_decomp} is a generalization of Corollary~\ref{idem_decomp_2} from idempotent elements to arbitrary non-nilpotent elements, a question arises whether Proposition~\ref{idem_decomp} could be generalized as well.
Do the summands in a minimal right decomposition of a non-nilpotent element form an orthogonal system? The answer is negative, and even more, there exist elements that possess no minimal right decomposition, whose summands would form an orthogonal system.
For example, in the ring $M_2(\mathbb{C})$, the matrix $A=E_{11}+E_{12}+E_{22}$ possesses no such decomposition (it is easy to verify that if $A=X+Y$, where $XY=YX=0$, then either $X=0$ or $Y=0$).
And even if an element possesses such a decomposition, this does not guarantee that all decompositions will have the same property. For example, the summands in the minimal right decomposition $B=E_{11}+2E_{22}$ form an orthogonal system (of non-nilpotent elements), while the summands in the decomposition $B=(E_{11}+E_{12})+(2E_{22}-E_{12})$ do not.

We are now ready to prove one of our main results, which gives a partial description of rank of idempotents in terms of invertible elements. It involves a condition, which resembles unit-regularity. 

\begin{theorem}\label{idem_rank}
For an idempotent $e \in R$ and a nonnegative integer $n$ the following conditions are equivalent:
\begin{enumerate}\sep
\item\label{normal} the right rank of $e$ is $n$,
\item\label{invertible} the right rank of $e$ is finite and greater than $n-1$, and for every $r \in R$ either the right rank of $er$ is less than $n$ or there exists $x \in \mathcal{U}(R)$, such that $er=ex$.
\end{enumerate}
\end{theorem}

\proof We first prove by induction that \ref{normal} implies \ref{invertible}. Since the right rank of $er$ is at most $n$, we only need to prove that if it is equal to $n$, then there exists an invertible element $x \in R$, such that $er=ex$. For $n=0$, i.e. $e=0$, this is obvious, just take $x=1$. Let $n \geq 1$ and suppose \ref{normal} implies \ref{invertible} for all smaller $n$.
Let $e=e_1+e_2+\ldots+e_n$ be a minimal right decomposition of $e$. By Proposition~\ref{idem_decomp}, elements $e_i$ are orthogonal idempotents. Let $f=e_2+e_3+\ldots+e_n$. By Corollary~\ref{idem_decomp_2}, $f$ is an idempotent of right rank $n-1$. Now suppose the right rank of $er$ is $n$ for some $r \in R$. Then $er=e_1r+fr$ implies that the right rank of $fr$ is $n-1$. By induction, there exists an invertible element $x \in R$, such that $fr=fx$. Hence
\begin{equation}\label{before}
erx^{-1}=(e_1r+fr)x^{-1}=e_1rx^{-1}+f.
\end{equation}

Assume first that $e_1rx^{-1}e_1 \neq 0$. Since $e_1$ is an idempotent of right rank $1$, $e_1Re_1$ is a division ring by Proposition~\ref{min_idem}. Hence, there exists $s \in R$, such that $e_1rx^{-1}e_1se_1=e_1se_1rx^{-1}e_1=e_1$. Observe that
\[e\big(e_1rx^{-1}+(1-e_1)\big)x=(e_1rx^{-1}+f)x=e_1r+fx=e_1r+fr=er,\]
thus it suffices to prove that the element $e_1rx^{-1}+(1-e_1)$ is invertible in $R$. We have
\begin{eqnarray*}
\big(e_1rx^{-1}+(1-e_1)\big)\big(e_1se_1+(1-e_1)\big) &=& e_1+e_1rx^{-1}(1-e_1)+(1-e_1)= \\
&=& 1+e_1rx^{-1}(1-e_1).
\end{eqnarray*}
Since $e_1se_1+(1-e_1)$ is invertible in $R$ with inverse $e_1rx^{-1}e_1+(1-e_1)$ and $1+e_1rx^{-1}(1-e_1)$ is invertible in $R$ with inverse $1-e_1rx^{-1}(1-e_1)$, the element $e_1rx^{-1}+(1-e_1)$ is invertible in $R$ as well.

Now assume that $e_1rx^{-1}e_1=0$. We prove by contradiction that in this case $e_1rx^{-1}(1-e) \neq 0$. Suppose $e_1rx^{-1}(1-e)=0$. Then
\[e_1rx^{-1}=e_1rx^{-1}(e+(1-e))=e_1rx^{-1}e=e_1rx^{-1}(e_1+f)=e_1rx^{-1}f.\]
Together with \eqref{before} this implies
\[erx^{-1}=e_1rx^{-1}+f=e_1rx^{-1}f+f=(1+e_1rx^{-1})f.\]
Hence, $\rank_r er=\rank_r erx^{-1} \leq \rank_r f=n-1$, which is a contradiction. Therefore $e_1rx^{-1}(1-e) \neq 0$. Since $e_1$ has right rank $1$, there exists $t \in R$, such that $e_1rx^{-1}(1-e)t=e_1$. Observe that
\[e\big(e_1rx^{-1}-(1-e)te_1+(1-e_1)\big)x=(e_1rx^{-1}+f)x=e_1r+fx=e_1r+fr=er,\]
thus it suffices to prove that the element $e_1rx^{-1}-(1-e)te_1+(1-e_1)$ is invertible in $R$. Since $e_1rx^{-1}(1-e)t=e_1$ and $e_1(1-e)=0$, which also implies $(1-e_1)(1-e)=1-e$, we have
\begin{eqnarray*}
&& \big(e_1rx^{-1}-(1-e)te_1+(1-e_1)\big)\big(1+(1-e)te_1\big)= \\
&& \hspace{3cm} =e_1rx^{-1}+e_1-(1-e)te_1+(1-e_1)+(1-e)te_1= \\
&& \hspace{3cm} =1+e_1rx^{-1}=1+e_1rx^{-1}(1-e_1),
\end{eqnarray*}
where the last equality follows from $e_1rx^{-1}e_1=0$. Since $1+(1-e)te_1$ is invertible in $R$ with inverse $1-(1-e)te_1$ and $1+e_1rx^{-1}(1-e_1)$ is invertible in $R$ with inverse $1-e_1rx^{-1}(1-e_1)$, the element $e_1rx^{-1}-(1-e)te_1+(1-e_1)$ is invertible in $R$ as well.

Now suppose \ref{invertible} holds. Then the right rank of $e$ is $m\geq n$. Let $e=e_1+e_2+\ldots+e_m$ be a minimal right decomposition of $e$, where $e_i$ are pairwise orthogonal idempotents by Proposition~\ref{idem_decomp}.
Since the right rank of $e(e-e_m)=e-e_m=e_1+e_2+\ldots+e_{m-1}$ is $m-1$ by Corollary~\ref{idem_decomp_2}, there is no invertible element $x \in R$ such that $e(e-e_m)=ex$, because the right rank of $ex$ is $m$ for every invertible element $x$. Hence \ref{invertible} implies that the right rank of $e(e-e_m)$ is less than $n$. Therefore $m<n+1$ and consequently $m=n$ as desired. \endproof

\begin{corollary}\label{idem_rank_c}
Let $e \in R$ be an idempotent of finite right rank $n$ and $r \in R$ an arbitrary element. If $er$ has right rank $n$, then there exists $x \in \mathcal{U}(R)$, such that $er=ex$.
\end{corollary}

\begin{remark}\label{rank_cond}
Theorem~\ref{idem_rank} does not yet give a complete characterization of the right rank of idempotents because of the presence of the condition that the right rank of $e$ is finite in item \ref{invertible} of the theorem.
This condition is redundant when $n=1$ but essential when $n \geq 2$. If $n=1$, the rest of item \ref{invertible} automatically implies that $eR$ is a minimal right ideal of $R$. Indeed, if $0 \neq a \in eR$, then there exist an invertible $x \in R$, such that $a=ea=ex$. The invertibility of $x$ implies $aR=eR$, so $eR$ is a minimal right ideal.
Now let $n\geq 2$ and consider the ring $R=T_2(\mathbb{C})$. Then the right rank of $E_{11}$ is infinite (see Example~\ref{mn_rank}). However, for every $A \in R$ either $E_{11}AE_{11}=0$, so that $E_{11}A$ has right rank at most $1$, or $E_{11}A=E_{11}(E_{11}A+E_{22})$, where $E_{11}A+E_{22}$ is invertible in $R$.
Nevertheless, we will be able to remove the aforementioned condition in case $R$ is a semiprime ring, thus producing a complete characterization of right rank (see Corollary~\ref{rank_char} below).
\end{remark}

As an application of the above results, we have the following nice proposition.

\begin{proposition}\label{unit_regular_prop}
Every regular element of finite right rank is unit-regular.
\end{proposition}

\proof Let $a$ be a regular element. Then there exists an element $b$ such that $a=aba$. Proposition~\ref{rank_ineq} implies that the idempotent $e=ab$ has the same right rank as $a$. By Corollary~\ref{idem_rank_c}, there exists an invertible element $x$ such that $a=ea=ex$, so $a$ is unit-regular. \endproof

By \cite[Remark~3.2]{Tug} every element in the socle of a semiprime ring is regular, hence Proposition~\ref{unit_regular_prop} immediately implies the following.

\begin{theorem}\label{unit_regular}
Every element in the socle of a unital semiprime ring is unit-regular.
\end{theorem}

Theorem~\ref{unit_regular} fails if $R$ is not semiprime, even for elements in the intersection of the left and right socle  of $R$. For example, consider again the ring $R=T_2(\mathbb{C})$. Then the matrix $E_{12}$ is contained in the intersection of left and right socle of $R$, however, $E_{12}$ also lies in the Jacobson radical of $R$, so it cannot be regular, let alone unit-regular, since the Jacobson radical does not contain nonzero idempotents.

Recall that the left and right rank in a semiprime ring coincide, so we omit the adjectives and simply speak of rank in this case. The next corollary gives a new inductive characterization of rank in semiprime rings in the sense of Theorem~\ref{idem_rank}. In particular, the condition discussed in Remark~\ref{rank_cond} is redundant here.

\begin{corollary}\label{rank_char}
Let $R$ be a semiprime ring and $n$ a nonnegative integer. For an element $a \in R$ the following conditions are equivalent:
\begin{enumerate}\sep
\item\label{normal_c} the rank of $a$ is $n$,
\item\label{invertible_c} the rank of $a$ is greater than $n-1$, and for every $r \in R$ either the rank of $ar$ is less than $n$ or there exists $x \in \mathcal{U}(R)$ such that $ar=ax$.
\end{enumerate}
\end{corollary}

\proof Suppose \ref{normal_c} holds. Then, by Theorem~\ref{unit_regular}, $a=ex$ for some idempotent $e$ and some invertible element $x$. Clearly, the rank of $e$ is $n$ as well. Take an arbitrary $r \in R$ and suppose the rank of $ar$ is $n$ (it cannot be greater than $n$). By Corollary~\ref{idem_rank_c}, there exists an invertible element $y$, such that $ar=e(xr)=ey=(ex)(x^{-1}y)=a(x^{-1}y)$. This proves \ref{invertible_c}, since $x^{-1}y$ is invertible.

Now suppose \ref{invertible_c} holds. We first prove that the rank of $a$ is finite. Observe that if $b \in aR$ is an element, such that $bR \varsubsetneq aR$, then the rank of $b$ is finite. Indeed, if the rank of $b=ar \in aR$ was infinite, condition \ref{invertible_c} would imply $b=ar=ax$ for some $x \in \mathcal{U}(R)$, but this would lead to a contradiction $aR=bx^{-1}R \subseteq bR$. To prove that the rank of $a$ is finite, we may assume that $\rank a \geq 2$, so that $aR$ is neither $0$ nor a minimal right ideal of $R$. Hence, there exists a right ideal $K \subseteq R$ such that $0 \neq K \varsubsetneq aR$. Choose a nonzero element $k \in K$. By the above, the rank of $k$ is finite. By Theorem~\ref{unit_regular}, $k=fy$ for some nonzero idempotent $f$ and some invertible element $y$. Hence $f=ky^{-1} \in K \subseteq aR$. Observe that $a=(1-f)a+fa$. The rank of $f$ is equal to the rank of $k$, so the element $fa$ has finite rank. Condition $f \in aR$ implies $(1-f)a \in aR$. If $(1-f)aR=aR$, then $f=f^2 \in f(aR)=f(1-f)aR=0$, which is a contradiction. Hence $(1-f)aR \varsubsetneq aR$ and consequently $(1-f)a$ has finite rank by the above. We conclude that $a=(1-f)a+fa$ has finite rank. By Theorem~\ref{unit_regular}, $a=ez$ for some idempotent $e$ and some invertible element $z$. The rank of $e$ is equal to the rank of $a$. It is easy to verify that the idempotent $e$ also satisfies condition \ref{invertible_c}, and since it has finite rank, it has rank $n$ by Theorem~\ref{idem_rank}. Hence $a$ has rank $n$ as well. \endproof

We believe that the characterization of rank in Corollary~\ref{rank_char} is new even in the case of square matrices over a field.

Here are some additional characterizations of rank in semiprime rings, which generalize some of the results obtained by Aupetit and Mouton \cite{Aup-Mou} and Bre\v sar and \v Semrl \cite{Bre-Sem} to semiprime rings.

\begin{proposition}
Let $R$ be a semiprime ring and $n$ a nonnegative integer. For an element $a \in R$ the following conditions are equivalent:
\begin{enumerate}
\item\label{char1} the rank of $a$ is $n$,
\item\label{char2} the right ideal $aR$ is a sum of $n$ minimal right ideals, but is not a sum of less than $n$ minimal right ideals,
\item\label{char3} the right ideal $aR$ is contained in $\soc R$ and has finite length $n$ as a right $R$-module.
\end{enumerate}
\end{proposition}

\proof Suppose $a$ has finite rank $n$. Clearly $aR$ cannot be a sum of less than $n$ minimal right ideals. By Theorem~\ref{unit_regular}, $a=eu$, where $u \in \mathcal{U}(R)$ and $e$ is an idempotent of rank $n$. Proposition~\ref{idem_decomp} implies $e=e_1+e_2+\ldots+e_n$, where $e_i$ are orthogonal idempotents of rank $1$. Since $e_i=ee_i \in eR$, we clearly have $aR=eR=e_1R+e_2R+\ldots+e_nR$, which implies \ref{char2}.

Suppose \ref{char2} holds and $aR=K_1+K_2+\ldots+K_n$, where $K_i$ are minimal right ideals of $R$. Then
\[0 \subset K_1 \subset K_1+K_2 \subset\ldots\subset K_1+K_2+\ldots+K_n\]
is a composition series of $eR$ (the second part of \ref{char2} implies that the composition factors are nonzero), so $aR$ has finite length $n$.

Now suppose \ref{char3} holds. By assumption $a$ has finite rank, say equal to $m$. From what we proved above it follows that the right ideal $aR$ has finite length $m$. By Jordan-H\" older Theorem the length is unique, hence $m=n$. \endproof

One interesting application of the results obtained above can be found in a subsequent paper \cite{Sto2}, where we discuss the structure of generalized corner rings.


\begin{thebibliography}{99}
\bibitem{Ami} S.A. Amitsur: Generalized polynomial identities and pivotal monomials, \emph{Trans. Amer. Math. Soc.} 114 (1965) 210--226. 
\bibitem{Aup-Mou} B. Aupetit, H. du T. Mouton: Trace and determinant in Banach algebras, \emph{Studia Math.} 121 (1996), no. 2, 115--136.
\bibitem{Bei-Mar-Mik} K.I. Beidar, W.S. Martindale  III, A.V. Mikhalev: \emph{Rings with generalized identities}, Marcel Dekker, Inc., New York, 1996.
\bibitem{Bre} M. Bre\v sar: \emph{Introduction to noncommutative algebra}, Springer International Publishing, Switzerland, 2014.
\bibitem{Bre-Sem} M. Bre\v sar, P. \v Semrl: Finite rank elements in semisimple Banach algebras, \emph{Studia Math.} 128 (1998), no. 3, 287--298.
\bibitem{Bri-Lin-Rau} R.M. Brits, L. Lindeboom, H. Raubenheimer: Rank and the Drazin inverse in Banach algebras, \emph{Studia Math.} 177 (2006), no. 3, 211--224.
\bibitem{Hun} T.W. Hungerford: \emph{Algebra}, Springer-Verlag, New York, 1974.
\bibitem{Lam} T.Y. Lam: \emph{A First Course in Noncommutative Rings, Second Edition}, Springer Science + Business Media, New York, 2001.
\bibitem{Mou-Rau} T. Mouton, H. Raubenheimer: On rank one and finite elements of Banach algebras, \emph{Studia Math.} 104 (1993), no. 3, 211--219.
\bibitem{Puh} J. Puhl: The trace of finite and nuclear elements in Banach algebras, \emph{Czechoslovak Math. J.} 28 (1978), no. 4, 656--676.
\bibitem{Row} L.H. Rowen: \emph{Polynomial Identities in Ring Theory}, Academic Press, Inc., New York, 1980.
\bibitem{Sto2} N. Stopar: Structure theorem for generalized corner rings, in preparation.
\bibitem{Tug} A. Tuganbaev: \emph{Rings Close to Regular}, Kluwer Academic Publishers, Dordrecht, 2002.
\end{thebibliography}
\end{document}